\def\N{\mathbb N}
\def\Z{\mathbb Z}

\def\Q{\mathbb Q}
\def\A{\mathcal A}

\def\pfz{\begin{proof}}
\def\pfk{\end{proof}}

\documentclass[runningheads,a4paper]{llncs}
\usepackage{amssymb,amsmath}

\usepackage{url}
\newcommand{\keywords}[1]{\par\addvspace\baselineskip
\noindent\keywordname\enspace\ignorespaces#1}

\newtheorem{lem}{Lemma}
\newtheorem{thm}[lem]{Theorem}

\newtheorem{coro}[lem]{Corollary}

\newtheorem{pozn}[lem]{Remark}

\begin{document}


\mainmatter  

\title{Ito-Sadahiro numbers vs.\ Parry numbers}


%
%
\author{Zuzana Mas\'akov\'a%
\and Edita Pelantov\'a}
%

\institute{Department of Mathematics FNSPE, Czech Technical University in Prague\\
Trojanova 13, 120 00 Praha 2, Czech Republic\\
zuzana.masakova@fjfi.cvut.cz, edita.pelantova@fjfi.cvut.cz}

%
%

\toctitle{Ito-Sadahiro numbers vs.\ Parry numbers}
\tocauthor{Z. Mas\'akov\'a, E. Pelantov\'a}
\maketitle

\begin{abstract}
We consider positional numeration system with negative base, as introduced by Ito and Sadahiro. In particular, we focus on algebraic properties of negative bases $-\beta$ for which the corresponding dynamical system is sofic, which happens, according to Ito and Sadahiro, if and only if the $(-\beta)$-expansion of $-\frac{\beta}{\beta+1}$ is eventually periodic. We call such numbers $\beta$ Ito-Sadahiro numbers and we compare their properties with Parry numbers, occurring in the same context for R\'enyi positive base numeration system.
\end{abstract}




\section{Introduction}

The expansion of a real number in the positional number system with base $\beta>1$, as defined by R\'enyi~\cite{Renyi} is closely related with the
transformation $T:[0,1)\mapsto[0,1)$, given by the prescription
$T(x):=\beta x-\lfloor\beta x\rfloor$. Every $x\in[0,1)$ is a sum of the infinite series
\begin{equation}\label{eq:1}
x=\sum_{i=1}^\infty\frac{x_i}{\beta^i}\,,\qquad\hbox{where }\ x_i=\lfloor\beta T^{i-1}(x)\rfloor\quad\hbox{for }i=1,2,3,\dots
\end{equation}
Directly from the definition of the transformation $T$ we can derive that the `digits' $x_i$ take values in the set $\{0,1,2,\cdots,\lceil\beta\rceil -1\}$ for $i=1,2,3,\cdots$.
The expression of $x$ in the form~\eqref{eq:1} is called the $\beta$-expansion of $x$. The number $x$ is thus represented by the infinite word $d_\beta(x)=x_1x_2x_3\cdots\in{\mathcal A}^\N$ over the alphabet $\A=\{0,1,2,\dots,\lceil\beta\rceil-1\}$.

From the definition of the transformation $\beta$ we can derive another important property, namely that the ordering on real numbers is carried over to
the ordering of $\beta$-expansions. In particular, we have for $x,y\in[0,1)$ that
$$
x\leq y \quad\iff\quad d_\beta(x) \preceq d_\beta(y)\,,
$$
where $\preceq$ is the lexicographical order on $\A^{\N}$, (ordering on the alphabet $\A$ is usual, $0<1<2<\cdots<\lceil\beta\rceil-1$).

In~\cite{Parry}, Parry has provided a criteria which decides whether an infinite word in $\A^\N$ is or not a $\beta$-expansion of some real number $x$.
The criteria is formulated using the so-called infinite expansion of $1$, denoted by $d^*_\beta(1)$, defined as a limit in the space $\A^\N$ equipped with the product topology, by
$$
d_\beta^*(1):=\lim_{\varepsilon\to0+}d_\beta(1-\varepsilon)\,.
$$
According to Parry, the string $x_1x_2x_3\cdots \in\A^\N$ represents the $\beta$-expansion of a number $x\in[0,1)$ if and only if
\begin{equation}\label{eq:2}
x_ix_{i+1}x_{i+2}\cdots \prec d_\beta^*(1)\quad \hbox{ for every }\ i=1,2,3,\cdots\,.
\end{equation}

The condition~\eqref{eq:2} ensures that the set ${\mathcal D}_\beta=\{d_\beta(x)\mid x\in[0,1)\}$ is shift invariant, and so the closure of ${\mathcal D}_\beta$ in $\A^\N$, denoted by $S_\beta$, is the subshift of the full shift $\A^\N$.

The notion of $\beta$-expansion can be naturally extended to all non-negative real numbers: The expression of a real number $y$ in the form
\begin{equation}\label{eq:3}
\begin{aligned}
y=y_k\beta^{k}+y_{k-1}\beta^{k-1}+y_{k-2}\beta^{k-2}+\cdots\,,\\[2mm]
\hbox{ where $k\in\Z$ and $y_ky_{k-1}y_{k-2}\cdots \in {\mathcal D}_\beta,$}
\end{aligned}
\end{equation}
is called the $\beta$-expansion of $y$.

Real numbers $y$ having in their $\beta$-expansion vanishing digits $y_i$ for all $i<0$ are usually called $\beta$-integers and the set of $\beta$-integers
is denoted by $\Z_\beta$. The notion of $\beta$-integers was first considered in~\cite{BuFrGaKr} as an aperiodic structure modeling non-crystallographic materials with long range order, called quasicrystals.

The choice of the base $\beta>1$ strongly influences the properties of $\beta$-expansions. It turns out that important role among bases  play such numbers $\beta$ for which $d_\beta^*(1)$ is eventually periodic. Parry, himself, has called these bases beta-numbers; now these numbers are commonly called Parry numbers. One can demonstrate the exceptional properties of Parry numbers on two facts:

\begin{itemize}
\item The subshift $S_\beta$ is sofic if and only if $\beta$ is a Parry number~\cite{ItoTakahashi}.

\item Distances between consecutive $\beta$-integers take finitely many values if and only if  $\beta$ is a Parry number~\cite{Thurston}.
\end{itemize}

Recently, Ito and Sadahiro~\cite{ItoSadahiro} suggested to study positional systems with negative base $-\beta$, where $\beta>1$. Representation of real numbers in such a system is defined using the transformation $T:[l_\beta,r_\beta)\mapsto[l_\beta,r_\beta)$, where $l_\beta=-\frac{\beta}{\beta+1}$, $r_\beta=1+l_\beta=\frac1{1+\beta}$. The transformation $T$ is defined by
\begin{equation}\label{eq:7}
T(x):=-\beta x -\lfloor-\beta x - l_\beta\rfloor\,.
\end{equation}
Every real $x\in I_\beta := [l_\beta,r_\beta)$ can be written as
\begin{equation}\label{eq:4}
x=\sum_{i=1}^\infty\frac{x_i}{(-\beta)^i}\,,\qquad\hbox{where }\ x_i=\lfloor-\beta T^{i-1}(x)-l_\beta\rfloor\quad\hbox{for }i=1,2,3,\dots
\end{equation}
The above expression is called the $(-\beta)$-expansion of $x$. It can also be written as the infinite word $d_{-\beta}(x)=x_1x_2x_3\cdots$.
One can easily show from~\eqref{eq:7} that the digits $x_i$, $i\geq 1$, take values in the set $\A=\{0,1,2,\dots,\lfloor\beta\rfloor\}$.
In this case, the ordering on the set of infinite words over the alphabet $\A$ which would correspond to the ordering of real numbers is the so-called alternate ordering: We say that $x_1x_2x_3\cdots \prec_{\hbox{\tiny alt}} y_1y_2y_3\cdots$ if for the minimal index $j$ such that $x_j\neq y_j$
it holds that $x_j(-1)^j<y_j(-1)^j$. In this notation, we can write for arbitrary $x,y\in I_\beta$ that
$$
x\leq y \quad\iff\quad d_{-\beta}(x) \preceq_{\hbox{\tiny alt}} d_{-\beta}(y)\,.
$$

In their paper, Ito and Sadahiro has provided a criteria to decide whether an infinite word $\A^\N$ belongs to the set of $(-\beta)$-expansions, i.e.\ to the set ${\mathcal D}_{-\beta}=\{d_{-\beta}(x) \mid x\in I_\beta\}$. This time, the criteria is given in terms of two infinite words, namely
$$
d_{-\beta}(l_\beta)\quad\hbox{ and }\quad d^*_{-\beta}(r_\beta):=\lim_{\varepsilon\to0+} d_{-\beta}(r_\beta-\varepsilon)\,.
$$
These two infinite words have close relation: If $d_{-\beta}(l_\beta)$ is purely periodic with odd period length, i.e.
$d_{-\beta}(l_\beta)=(d_1d_2\cdots d_{2k+1})^\omega$, then we put $d_{-\beta}^*(r_\beta)=\big(0d_1d_2\cdots (d_{2k+1}-1)\big)^\omega$.
(As usual, the notation $w^\omega$ stands for infinite repetition of the string $w$.)
In all other cases one has $d_{-\beta}^*(r_\beta)=0d_{-\beta}(l_\beta)$.

Ito and Sadahiro have shown that an infinite word $x_1x_2x_3\cdots$ represents a $(-\beta)$-expansion of some $x\in[l_\beta,r_\beta)$ if and only if
for every $i\geq 1$ it holds that
\begin{equation}\label{eq:5}
d_{-\beta}(l_\beta) \preceq_{\hbox{\tiny alt}} x_ix_{i+1}x_{i+2}\cdots \prec_{\hbox{\tiny alt}} d_{-\beta}^*(r_\beta)\,.
\end{equation}
The above condition ensures that the set ${\mathcal D}_{-\beta}$ of infinite words representing $(-\beta)$-expansions is shift invariant.
In~\cite{ItoSadahiro} it is shown that the closure of ${\mathcal D}_{-\beta}$ defines a sofic system if and only if $d_{-\beta}(l_\beta)$ is
eventually periodic.

In analogy to the definition of Parry numbers, we suggest that numbers $\beta>1$ such that $d_{-\beta}(l_\beta)$ is eventually periodic be called
Ito-Sadahiro numbers.
The relation of the set of Ito-Sadahiro numbers and the set of Parry numbers is not obvious. We do not know any example of a Parry number which is not an Ito-Sadahiro number or vice-versa. Bassino~\cite{Bassino} has shown that quadratic and cubic not totally real numbers $\beta$ are Parry numbers if and only if they are Pisot numbers. For the same class of numbers, we prove in~\cite{MaPeVa} that $\beta$ is Ito-Sadahiro if and only if it is Pisot. This means that notions of Parry numbers and Ito-Sadahiro numbers on the mentioned type of irrationals coincide.

In this paper we study numbers with eventually periodic $(-\beta)$-expansion. Statements which we show, as well as results of other authors we recall,
demonstrate similarities between the behaviour of $\beta$-expansions and $(-\beta)$-expansions. We mention also phenomena in which the two essentially differ. Nevertheless, we are in favour of the hypothesis that the set of Parry numbers and the set of Ito-Sadahiro numbers coincide.

\section{Preliminaries}

Let us first recall some number theoretical notions. A complex number $\beta$ is called an algebraic number, if it is root of a monic polynomial $x^n + a_{n-1}x^{n-1}+\cdots +a_1x + a_0$, with rational coefficients $a_0,\dots,a_{n-1}\in\Q$. Monic polynomial with rational coefficients and root $\beta$ of
the minimal degree among all polynomials with the same properties is called the minimal polynomial of $\beta$, its degree is called the degree of $\beta$. The roots of the minimal polynomial are algebraic conjugates.

If the minimal polynomial of $\beta$ has integer coefficients, $\beta$ is called an algebraic integer. An algebraic integer $\beta>1$ is called a Perron number, if all its conjugates are in modulus strictly smaller than $\beta$. An algebraic integer $\beta>1$ is called a Pisot number, if all its conjugates are in modulus strictly smaller than $1$. An algebraic integer $\beta>1$ is called a Salem number, if all its conjugates are in modulus smaller than or equal to $1$ and $\beta$ is not a Pisot number.

If $\beta$ is an algebraic number of degree $n$, then the minimal subfield of the field of complex numbers containing $\beta$ is denoted by $\Q(\beta)$
and is of the form
$$
\Q(\beta) = \{c_0+c_1\beta + \cdots + c_{n-1}\beta^{n-1} \mid c_i\in\Q\}\,.
$$
If $\gamma$ is a conjugate of an algebraic number $\beta$, then the fields $\Q(\beta)$ and $\Q(\gamma)$ are isomorphic. The corresponding isomorphism is given by the prescription
$$
c_0+c_1\beta + \cdots + c_{n-1}\beta^{n-1} \quad\mapsto\quad c_0+c_1\gamma + \cdots + c_{n-1}\gamma^{n-1}\,.
$$
In particular, it means that $\beta$ is a root of some polynomial $f$ with rational coefficients if and only if $\gamma$ is a root of the same polynomial $f$.

\section{Ito-Sadahiro polynomial}

From now on, we shall consider for bases of the numeration system only Ito-Sadahiro numbers, i.e.\ numbers $\beta$ such that
\begin{equation}\label{eq:6}
d_{-\beta}(l_\beta) = d_1\cdots d_m (d_{m+1}\cdots d_{m+p})^\omega\,.
\end{equation}
Without loss of generality we shall assume that $m\geq 0$, $p\geq 1$ are minimal values so that $d_{-\beta}(l_\beta)$ can be written in the above form.
Recall that $l_\beta=-\frac{\beta}{\beta+1}$. Therefore~\eqref{eq:6} can be rewritten as
$$
-\frac{\beta}{\beta+1} = \frac{d_1}{-\beta} + \cdots + \frac{d_m}{(-\beta)^m} +
\left(\frac{d_{m+1}}{(-\beta)^{m+1}} + \cdots + \frac{d_{m+p}}{(-\beta)^{m+p}}\right) \sum_{i=0}^\infty \frac1{(-\beta)^{p\,i}}\,,
$$
and after arrangement
$$
0 = \frac{-\beta}{-\beta-1} + \frac{d_1}{-\beta} + \cdots + \frac{d_m}{(-\beta)^m} + \frac{(-\beta)^p}{(-\beta)^p-1}
\left(\frac{d_{m+1}}{(-\beta)^{m+1}} + \cdots + \frac{d_{m+p}}{(-\beta)^{m+p}}\right) \,.
$$
Multiplying by $(-\beta)^m\big((-\beta)^p-1\big)$, we obtain the following lemma.

\begin{lem}
Let $\beta$ be an Ito-Sadahiro number and let $d_{-\beta}(l_\beta)$ be of the form~\eqref{eq:6}. Then $\beta$ is a root of the polynomial
\begin{equation}\label{eq:ISpolyObec}
\begin{aligned}
P(x) &= (-x)^{m+1}\sum_{i=0}^{p-1}(-x)^i + \big((-x)^p-1\big) \sum_{i=1}^m d_i (-x)^{m-i} + \sum_{i=m+1}^{m+p} d_i (-x)^{m+p-i}\,.
\end{aligned}
\end{equation}
Such polynomial is called the Ito-Sadahiro polynomial of $\beta$.
\end{lem}

\begin{coro}\label{c}
An Ito-Sadahiro number is an algebraic integer of degree smaller than or equal to $m+p$, where $m$, $p$ are given by~\eqref{eq:6}.
\end{coro}

It is useful to mention that the Ito-Sadahiro polynomial is not necessarily irreducible over $\Q$. As an example one can take the minimal Pisot number. For such $\beta$, we have $d_{-\beta}(l_\beta)=1001^\omega$, and thus the Ito-Sadahiro polynomial is equal to $P(x)=x^4-x^3-x^2+1=(x-1)(x^3-x-1)$, where
$x^3-x-1$ is the minimal polynomial of $\beta$.

\begin{pozn}
Note that for $p=1$ and $d_{m+1}=0$, we have $d_{-\beta}(l_\beta)=d_1\cdots d_m0^\omega$, and the Ito-Sadahiro polynomial of $\beta$ is of the form
\begin{equation}\label{eq:ISpolySimple}
P(x) = (-x)^{m+1} + d_1(-x)^m + (d_2-d_1)(-x)^{m-1} + \cdots + (d_m-d_{m-1})(-x) - d_m\,,
\end{equation}
and thus $\beta$ is an algebraic integer of degree at most $m+1$.
\end{pozn}

\begin{thm}\label{thm:1}
Let $\beta$ be an Ito-Sadahiro number. All roots $\gamma$, $\gamma\neq\beta$, of the Ito-Sadahiro polynomial (in particular all conjugates of $\beta$) satisfy $|\gamma|< 2$.
\end{thm}

\pfz
Since $\beta$ is a root of the Ito-Sadahiro polynomial $P$, there must exist a polynomial $Q$ such that
$P(x)=(x-\beta)Q(x)$. Let us first determine $Q$ and show that it is a monic polynomial with coefficients in modulus not exceeding 1.
The coefficients $d_i$ in the polynomial $P$ in the form~\eqref{eq:ISpolyObec} are the digits of the $(-\beta)$-expansion of $l_\beta$, and thus,
by~\eqref{eq:4}, they satisfy $d_i=\lfloor-\beta T^{i-1}(l_\beta)-l_\beta\rfloor$.
Relation~\eqref{eq:7} then implies $T^{i}(l_\beta)=-\beta T^{i-1}(l_\beta) - \lfloor-\beta T^{i-1}(l_\beta)-l_\beta\rfloor$, wherefrom we have
$$
d_i=- T^{i}(l_\beta) -\beta T^{i-1}(l_\beta)\,.
$$
For simplicity of notation in this proof, denote $T_i=T^{i}(l_\beta)$, for $i=0,1,\dots,m+p$. Substituting $d_i=-T_i-\beta T_{i-1}$ into~\eqref{eq:ISpolyObec},
we obtain
\begin{equation}\label{eq:8}
\begin{aligned}
P(x) &= (-x)^{m+1}\sum_{i=0}^{p-1}(-x)^i + \big((-x)^p-1\big) \sum_{i=1}^m (-T_i-\beta T_{i-1}) (-x)^{m-i} + \\
&\hspace*{3.3cm} + \sum_{i=m+1}^{m+p} (-T_i-\beta T_{i-1}) (-x)^{m+p-i} = \\
&= (-x)^{m+1}\sum_{i=0}^{p-1}(-x)^i + \big((-x)^p-1\big) (x-\beta)\sum_{i=2}^m T_{i-1} (-x)^{m-i} + \\
&\hspace*{0.5cm} + (x-\beta)\sum_{i=1}^{p} T_{m+i-1} (-x)^{p-i}
- \big((-x)^p-1\big) \beta T_0 (-x)^{m-1} + \\
&\hspace*{5cm}+ \ T_m - T_{m+p}\,.
\end{aligned}
\end{equation}
First realize that $T_m - T_{m+p}=0$, since $d_{-\beta}(l_\beta)$ is eventually periodic with preperiod of length $m$ and period of length $p$.
As $T_0=T^0(l_\beta)=-\frac{\beta}{\beta+1}$, we can derive that
$$
(-x)^{m+1}\sum_{i=0}^{p-1}(-x)^i - \big((-x)^p-1\big) \beta T_0 (-x)^{m-1} = (-x)^{m-1}(x-\beta)(x-T_0) \sum_{i=0}^{p-1}(-x)^i\,.
$$
Putting back to~\eqref{eq:8}, we obtain that the desired polynomial $Q$ defined by $P(x)=(x-\beta)Q(x)$ is of the form
$$
\begin{aligned}
Q(x) &= (-x)^{m-1}(x-T_0) \sum_{i=0}^{p-1}(-x)^i + \\
&\hspace*{0.4cm} + \big((-x)^p-1\big)\sum_{i=2}^m T_{i-1} (-x)^{m-i} + \sum_{i=1}^{p} T_{m+i-1} (-x)^{p-i}\,,
\end{aligned}
$$
which can be rewritten in another form, namely,
\begin{equation}\label{eq:Q}
\begin{aligned}
Q(x)
&=  -(-x)^{m+p-1} + \sum_{i=m}^{m+p-2}(T_{m+p-1-i}-T_0-1)(-x)^i + \\
&\hspace*{2.5cm}+ \sum_{i=0}^{m-1} (T_{m+p-1-i}-T_{m-1-i}) (-x)^i\,.
\end{aligned}
\end{equation}
Note that the coefficients at individual powers of $-x$ are of two types,
namely
$$
T_{m+p-1-i}-T_0-1 \in [-1,0)\,,\quad\hbox{ and }\quad T_{m+p-1-i}-T_{m-1-i} \in (-1,1)\,.
$$
In order to complete the proof, realize that every root $\gamma$, $\gamma\neq\beta$, of the polynomial $P$ satisfies $Q(\gamma)=0$.
We thus have
$$
(-\gamma)^{m+p-1} =  \sum_{i=m}^{m+p-2}(T_{m+p-1-i}-T_0-1)(-\gamma)^i + \sum_{i=0}^{m-1} (T_{m+p-1-i}-T_{m-1-i}) (-\gamma)^i\,,
$$
and hence
$$
|\gamma|^{m+p-1} \leq 
\sum_{i=0}^{m+p-2}|\gamma|^i = \frac{|\gamma|^{m+p-1}-1}{|\gamma|-1} < \frac{|\gamma|^{m+p-1}}{|\gamma|-1}\,.
$$
From this, one easily derives that $|\gamma|<2$.
\pfk

\begin{coro}\label{cc}
Every Ito-Sadahiro number $\beta\geq 2$ is a Perron number.
\end{coro}

%
%

\section{Periodic expansions in the Ito-Sadahiro system}

Representations of numbers in the numeration system with negative base from the point of view of dynamical systems has been studied by Frougny and Lai~\cite{ChiaraFrougny}. They have shown the following statement.

\begin{thm}\label{thm:FL}
If $\beta$ is a Pisot number, then $d_{-\beta}(x)$ is eventually periodic for any $x\in I_\beta\cap\Q(\beta)$.
\end{thm}

In particular, their result implies that every Pisot number is an Ito-Sadahiro number. Here, we show a `reversed' statement.

\begin{thm}\label{thm:nas}
If any $x\in I_\beta\cap\Q(\beta)$ has eventually periodic $(-\beta)$-expansion, then $\beta$ is either Pisot or Salem
number.
\end{thm}

\pfz
First realize that since $l_{-\beta}\in\Q(\beta)$, by assumption, $d_{-\beta}(l_\beta)$ is eventually periodic, and thus $\beta$ is an Ito-Sadahiro number. Therefore using Corollary~\ref{c}, $\beta$ is an algebraic integer. It remains to show that all conjugates of $\beta$ are in modulus smaller than or equal to $1$.

Consider a real number $x$ whose $(-\beta)$-expansion is of the form $d_{-\beta}(x)=x_1x_2x_3\cdots$. We now show that
\begin{equation}\label{eq:aux}
x_1=x_2=\cdots =x_{k-1}=0\quad\hbox{and}\quad x_k\neq0\quad\hbox{ implies }\quad
|x|\geq \frac1{\beta^k(\beta+1)}\,.
\end{equation}
In order to see this, we estimate the series
$$
|x|=\Big|\frac{x_k}{(-\beta)^k} + \sum_{i=1}^\infty\frac{x_{k+i}}{(-\beta)^{k+i}}\Big| \geq \frac1{\beta^k} - \frac1{\beta^k}\Big|\sum_{i=1}^\infty\frac{x_{k+i}}{(-\beta)^i}\Big|\,.
$$
Since the set ${\mathcal D}_{-\beta}$ of all $(-\beta)$-expansions is shift invariant, the sum $\sum_{i=1}^\infty\frac{x_{k+i}}{(-\beta)^i}$
is a $(-\beta)$-expansion of some $y\in I_\beta$. Therefore we can write
$$
|x|\geq \frac1{\beta^k} - \frac1{\beta^k}|y| \geq \frac1{\beta^k} - \frac1{\beta^k} \frac{\beta}{\beta+1} = \frac{1}{\beta^k(\beta+1)}\,.
$$
As $\beta>1$, there exists $L\in\N$ such that
$$
-\frac{\beta}{\beta+1} < \frac1{(-\beta)^{2L+1}}\,.
$$
Let $M\in\N$ satisfy $M>2L+1$. Choose a rational number $r$ such that
\begin{equation}\label{eq:pf4}
\frac1{(-\beta)^{2L+1}} < r < \frac1{(-\beta)^{2L+1}} + \frac{1}{\beta^M(\beta+1)}\,.
\end{equation}
According to the auxiliary statement~\eqref{eq:aux}, the $(-\beta)$-expansion of $r$ must be of the form
\begin{equation}\label{eq:pf1}
r=\frac1{(-\beta)^{2L+1}} + \sum_{i=M+1}^\infty \frac{r_i}{(-\beta)^i}\,.
\end{equation}
As $r$ is rational, by assumption, the infinite word $r_{M+1}r_{M+2}\cdots$ is eventually periodic and by summing a geometric series,
$\sum_{i=M+1}^\infty \frac{r_i}{(-\beta)^i}$ can be rewritten as
$$
\sum_{i=M+1}^\infty \frac{r_i}{(-\beta)^i} = c_0+c_1\beta+\cdots +c_{n-1}\beta^{n-1}\in\Q(\beta)\,,
$$
where $n$ is the degree of $\beta$.

In order to prove the theorem by contradiction, assume that a conjugate $\gamma\neq \beta$ is in modulus greater than 1. By application of the
isomorphism between $\Q(\beta)$ and $\Q(\gamma)$, we get
$$
c_0+c_1\gamma+\cdots +c_{n-1}\gamma^{n-1}  = \sum_{i=M+1}^\infty \frac{r_i}{(-\gamma)^i}\,,
$$
and thus
\begin{equation}\label{eq:pf2}
r= \frac1{(-\gamma)^{2L+1}} +\sum_{i=M+1}^\infty \frac{r_i}{(-\gamma)^i}\,.
\end{equation}
Subtracting~\eqref{eq:pf2} from~\eqref{eq:pf1}, we obtain
\begin{equation}\label{eq:pf3}
0< \Big|\frac{1}{(-\beta)^{2L+1}} - \frac1{(-\gamma)^{2L+1}}\Big| \leq \sum_{i=M+1}^\infty r_i \big|(-\beta)^{-i}-(-\gamma)^{-i}\big|\leq
2\lfloor\beta\rfloor \frac{\eta^{M+1}}{1-\eta}\,,
\end{equation}
where $\eta = \max\{|\beta|^{-1},|\gamma|^{-1}\}<1$.
Obviously, for any $M>2L+1$, we can find a rational $r$ satisfying~\eqref{eq:pf4} and thus derive the inequality~\eqref{eq:pf3}. However,
the left-hand side of~\eqref{eq:pf3} is a fixed positive number, whereas the right-hand side decreases to zero with increasing $M$,
which is a contradiction.
\pfk

In order to stress the analogy of the Ito-Sadahiro numeration system with R\'enyi $\beta$-expansions of numbers, recall that already
Schmidt in~\cite{schmidt} has shown that for a Pisot number $\beta$, any $x\in[0,1)\cap\Q(\beta)$ has an eventually periodic $\beta$-expansion
and also the converse, that every $x\in[0,1)\cap\Q(\beta)$ having an eventually periodic $\beta$-expansion forces $\beta$ to be either Pisot or Salem number.
In fact, the proof of Theorem~\ref{thm:FL} given by Frougny and Lai, as well as our proof of Theorem~\ref{thm:nas} are using the ideas presented in~\cite{schmidt}.

A special case of numbers with periodic $(-\beta)$-expansion is given by those numbers $x$ for which the infinite word $d_{-\beta}(x)$
has suffix $0^\omega$. We then say that the expansion $d_{-\beta}(x)$ is finite. An example of such a number is $x=0$ with $(-\beta)$-expansion $d_{-\beta}(x)=0^\omega$. As it is shown in~\cite{MaPeVa}, if $\beta<\frac12(1+\sqrt5)$, then $x=0$ is the only number with finite $(-\beta)$-expansion.
This property of the Ito-Sadahiro numeration system has no analogue in R\'enyi $\beta$-expansions; for positive base, the set of finite $\beta$-expansions
is always dense in $[0,1)$.

Just as in the numeration system with positive base, we can extend the definition of $(-\beta)$-expansions of $x$ to all real numbers $x$, and define
the notion of a $(-\beta)$-integer as a real number $y$ such that
$$
y=y_k(-\beta)^k + \cdots + y_1(-\beta) + y_0\,,
$$
where  $y_k\cdots y_1y_00^\omega$ is the $(-\beta)$-expansion of some number in $I_\beta$. The set of $(-\beta)$-integers is denote by $\Z_{-\beta}$. In this notation, we can write for the set of all numbers with finite $(-\beta)$-expansions
$$
{\rm Fin}(-\beta) = \bigcup_{k=0}^\infty\frac1{(-\beta)^k} \Z_{-\beta}\,.
$$
It is not surprising that arithmetical properties of $\beta$-expansions and $(-\beta)$-expan\-sions depend on the choice of the base $\beta$. It can be shown that both $\Z_\beta$ and $\Z_{-\beta}$ is closed under addition and multiplication if and only if $\beta\in\N$. On the other hand, ${\rm Fin}(\beta)$
and ${\rm Fin}(-\beta)$ can have a ring structure even if $\beta$ is not an integer. Frougny and Solomyak~\cite{FruSo} have shown that if ${\rm Fin}(\beta)$
is a ring, then $\beta$ is a Pisot number. Similar result is given in~\cite{MaPeVa} for negative base: ${\rm Fin}(-\beta)$
being a ring implies that $\beta$ is either a Pisot or a Salem number. In~\cite{MaPeVa} we also prove the conjecture of Ito and Sadahiro that in case of quadratic Pisot base $\beta$ the set ${\rm Fin}(-\beta)$ is a ring if and only if the conjugate of $\beta$ is negative.

\section{Comments and open questions}

\begin{itemize}
\item
Every Pisot number is a Parry number and every Parry number is a Perron number, and neither of the statements can be reversed. The former is a consequence of the mentioned result of Schmidt, the latter statement follows for example from the fact that every Perron number has an associated canonical substitution $\varphi_\beta$, see~\cite{Fabre}. The substitution is primitive, and its incidence matrix has $\beta$ as its eigenvalue. The fixed point of $\varphi_\beta$ is an infinite word which codes the sequence of distances between consecutive $\beta$-integers.

\item
For the negative base numberation system, we can derive from Theorem~\ref{thm:FL} that every Pisot number is an Ito-Sadahiro number. From Corollary~\ref{cc} we know  that an Ito-Sadahiro number $\beta\geq 2$ is a Perron number. Based on our investigation, we conjecture that for any
Ito-Sadahiro number $\beta\geq \frac12(1+\sqrt5)$, the sequence of distances between consecutive $(-\beta)$-integers can be coded by a fixed point of a `canonical' substitution which is primitive and its incidence matrix has $\beta^2$ for its dominant eigenvalue. Thus we expect that every Ito-Sadahiro number $\beta\geq \frac12(1+\sqrt5)$ is also a Perron number. In case that $\beta<\frac12(1+\sqrt5)$, we have $\Z_{-\beta}=\{0\}$ and so
the situation is not at all obvious.

\item
In~\cite{solomyak}, Solomyak has explicitely described the set of conjugates of all Parry numbers. In particular, he has shown that this set is included in the complex disc of radius $\frac12(1+\sqrt5)$, and that this radius cannot be diminished. For his proof it was important that all conjugates
of a Parry number are roots of a polynomial with real coefficients in the interval $[0,1)$. In the proof of Theorem~\ref{thm:1} we show that conjugates of an Ito-Sadahiro number are roots of a polynomial~\eqref{eq:Q} with coefficients in $[-1,1]$. From this, we derive that conjugates of Ito-Sadahiro numbers lie in the complex disc of radius $\leq 2$. We do not know whether this value can be diminished.

\end{itemize}

\subsubsection*{Acknowledgments.}

We acknowledge financial support by the Czech Science Foundation
grant 201/09/0584 and by the grants MSM6840770039 and LC06002 of
the Ministry of Education, Youth, and Sports of the Czech
Republic.


\end{document}